\documentclass[12pt,twoside]{article}
\usepackage{geometry}
\usepackage[utf8]{inputenc}
\usepackage[T1]{fontenc}
\usepackage[english]{babel}
\usepackage{amsmath, amssymb, amsthm}            
\usepackage{amstext, amsfonts, a4}
\usepackage{graphicx}
\usepackage{hyperref}
\usepackage{color}
\usepackage{stmaryrd}
\usepackage{todonotes}
\usepackage{multirow}

\theoremstyle{plain}
\newtheorem{Theorem}{Theorem}[section]

\newcommand{\RR}{\mathbb R}
\newcommand{\QQ}{\mathbb Q}

\DeclareMathOperator{\GL}{\operatorname{GL}}
\DeclareMathOperator{\Prym}{\operatorname{Prym}}

\title{There are no primitive Teichm\"uller curves in $\Prym(2,2)$}
\author{\textsc{Julien Boulanger and Sam Freedman}}
\begin{document}

\maketitle

\begin{abstract}
We complete the work of Lanneau--M\"oller \cite{LM18} to show that there are no primitive Teichm\"uller curves in $\Prym(2,2)$.
\end{abstract}

\section{Introduction}
Teichm\"uller curves are closed $\GL(2, \RR)$-orbits in the moduli space of translation surfaces $\Omega\mathcal{M}_g$ that descend to isometrically-immersed algebraic curves in the moduli space of genus $g$ Riemann surfaces $\mathcal{M}_g$.
Lanneau--M\"oller \cite{LM18} searched for \emph{geometrically primitive} Teichm\"uller curves, i.e., those not arising from a covering construction, in a certain locus $\Prym(2,2)$ of genus 3 translation surfaces with two cone-points called \emph{Prym eigenforms}.
In the spirit of McMullen's proof that the decagon is the unique primitive Teichm\"uller curve in $\Omega \mathcal{M}_2(1, 1)$  (see \cite[Theorem 6.3]{McMullen2006}), they reduced their search to considering whether any of 92 candidate Prym eigenforms generate a Teichm\"uller curve.
In this note, we analyze those 92 candidates and show:
\begin{Theorem}\label{thm:main_thm}
The Prym locus $\Prym(2,2) \subset \Omega \mathcal{M}_3$ contains no primitive Teichm\"uller curves.
\end{Theorem}
Our main tool is the newly-available SageMath package \emph{Flatsurf} which partially computes the $\GL(2, \RR)$-orbit closure of an input translation surface.
(See Appendix B of Delecroix--R\"uth--Wright \cite{DRW21} for an introduction to \emph{Flatsurf}.)
Using \emph{Flatsurf}, we construct the 92 candidates Prym eigenforms and find a periodic direction that violates the Veech dichotomy.
See \S\S 2.2 for details.

We thank Erwan Lanneau for his assistance with this project.

\subsection*{Background}
Classifying the translation surfaces that generate Teichm\"uller curves, called \emph{Veech surfaces}, is a central problem in Teichm\"uller dynamics.
See Hubert--Schmidt \cite{HubertSchmidt06} for an introduction to Veech surfaces, M\"oller \cite{MollerTCurves} for a list of known examples of Teichm\"uller curves, Lanneau--M\"oller \cite{LM18} for a history of classification results, and McMullen \cite{McMullen2021BilliardsAT} for a general survey.

While McMullen \cite{McMullen2007} classified all primitive Teichm\"uller curves in genus two, we do not have a complete classification of primitive Teichm\"uller curves in all other genera.
In genus three, Bainbridge--Habegger--M\"oller \cite{BHM} found a numerical bound on the number of \emph{algebraically primitive} Teichm\"uller curves, yet the bound is too large for a complete classification.
(An \emph{algebraically primitive} Teichm\"uller curve is one for which the trace field of the generating Veech surface has maximal degree, i.e., degree the genus of the Veech surface.
We remark that algebraically primitive Teichm\"uller curves are automatically geometrically primitive, yet the converse does not hold.)
In fact, all but finitely many geometrically primitive Teichm\"uller curves in genus 3 lie in the Prym loci (see McMullen \cite[Theorem 5.5]{McMullen2021BilliardsAT}).

At another level, Lanneau--M\"oller \cite{LM18} began searching for Teichm\"uller curves among translation surfaces having a certain \emph{Prym involution} that generalizes the hyperelliptic involution in genus 2.
(See Lanneau--M\"oller \cite{LM18} for the definition.)
They showed that there are no \emph{geometrically primitive}
Teichm\"uller curves in $\Prym(2,1,1)$ and initiated the proof for $\Prym(2,2)$ that we finish in this note.
The remaining case of whether there are primitive Teichm\"uller curves in the locus $\Prym(1,1,1,1)$ is still open and seems interesting. 
Outside the Prym loci, Winsor \cite{Winsor22} showed that the Veech 14-gon generates the unique algebraically primitive Teichm\"uller curve in the hyperelliptic component of the stratum $\Omega \mathcal{M}_3(2, 2)$.

\section{Proof of \autoref{thm:main_thm}.}
\subsection{Candidate surfaces}
In this section we review, following Lanneau--M\"oller \cite{LM18}, how to construct the 92 candidate surfaces that could generate a Teichm\"uller curve.

\begin{figure}[tb]
\centering
\includegraphics[scale=0.9]{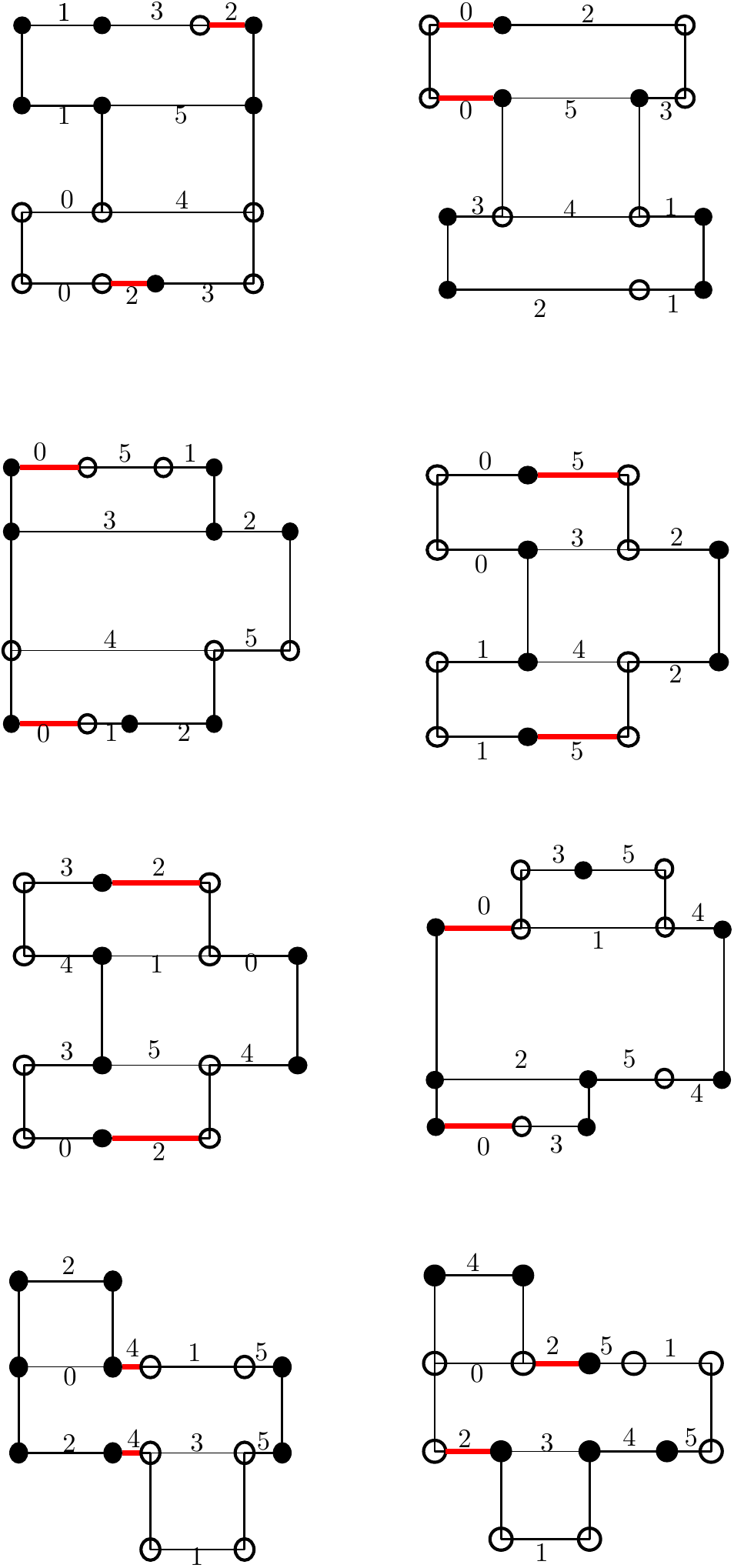}
\caption{The 8 polygonal models. The slit parameter encodes the length of the bold red relative period.}
\label{fig:models}
\end{figure}

Every surface in $\Prym(2,2)$ decomposes into horizontal cylinders (after possibly rotating the surface).
The combinatorics of the horizontal separatrices matches one of the eight polygonal models in \autoref{fig:models}.
(See also Figure 3 in Lanneau--M\"oller \cite{LM18}.)
In each underlying polygonal model, the Prym involution fixes one cylinder, called $C_2$, and permutes the other two, called $C_1$ and $C_3$.
The surface is then determined by choices of the widths, heights and twists of the cylinders $C_1$ and $C_2$, as well as of the length of a relative period called the \emph{slit parameter} $s$.
See Proposition 8.1 in Lanneau--M\"oller \cite{LM18} for an example of one of the eight possible separatrix diagrams.

In this respect, Lanneau--M\"oller \cite{LM18} determine a finite list of possible values for the width and the height of $C_2$, the surface being normalized with $w(C_1)=h(C_1)=1$.
More precisely, they first show (see their Theorem 4.5 and \S 6.2) that the only possible trace fields are $\QQ[\sqrt{D}]$ for $D \in \{2,3,33\}$.
They then encode the data of the widths and the heights of the horizontal and vertical cylinders in a \emph{reduced intersection matrix}.
In Proposition 5.5 and \S\S 6.3--6.5, they compute a finite list of such matrices that could give rise to Teichm\"uller curves; we recall this list in \autoref{table:table1} along with the associated width and height of $C_2$.

\begin{table}[tb]
\center
\begin{tabular}{|c|c|c|c|}
  \hline
 \multirow{2}{*}{Trace field} & Reduced & \multirow{2}{*}{$w(C_2)$} & \multirow{2}{*}{$h(C_2)$}\\
 & intersection matrix & & \\
\hline
  $\QQ[\sqrt{2}]$ &
  $\begin{pmatrix} 72 & 48\\ 24 & 18 \end{pmatrix}$ &
  $\frac{\sqrt{2}}{2}$ & $2\sqrt{2}$ \\
\hline
  $\QQ[\sqrt{3}]$ &
  $\begin{pmatrix} 72 & 24\\ 12 & 6 \end{pmatrix}$ &
  $\frac{-1+\sqrt{3}}{2}$ & $-2+2\sqrt{3}$ \\
\hline
  $\QQ[\sqrt{3}]$ & 
  $\begin{pmatrix} 72 & 24\\ 48 & 18 \end{pmatrix}$ &
  $\frac{1+\sqrt{3}}{2}$ & $2+2\sqrt{3}$ \\
\hline
  $\QQ[\sqrt{3}]$ &
  $\begin{pmatrix} 36 & 12\\ 30 &12 \end{pmatrix}$ & $\sqrt{3}$ & $\frac{2\sqrt{3}}{3}$ \\
\hline
  $\QQ[\sqrt{33}]$ &
  $\begin{pmatrix} 6 & 24\\ 12 & 54 \end{pmatrix}$ & $\frac{3+\sqrt{33}}{2}$ & $\frac{3+\sqrt{33}}{6}$\\
\hline
\end{tabular}
\caption[]{
The list of trace fields and reduced matrices that could give rise to Teichm\"uller curves\footnotemark, and the corresponding width and height of the cylinder $C_2$.
}
\label{table:table1}
\end{table}
\footnotetext{The last two reduced matrices of Lanneau--M\"oller\cite{LM18} do not generate Teichm\"uller curves according to Table 2 of \cite{LM18}, so we do not consider them.}

Having established the width and height parameters of the cylinders, it remains to control the length of the slit parameter and the twists of the cylinders.
In \S 8.2 of Lanneau--M\"oller \cite{LM18} and the corresponding code, they explicitly enumerate the finite list of possible values for such parameters that would lead to one of each admissible reduced intersection matrix.
This amounts to a finite list of candidate surfaces (see Table 2 of Lanneau--M\"oller \cite{LM18}).
In the code joint to this note, we list all those candidate surfaces and construct them using \emph{Flatsurf}.

\subsection{Cylinder decompositions}
Given a translation surface $M$, Flatsurf can efficiently generates many saddle connections on $M$.
For each such saddle connection, it can compute the decomposition of the straight-line flow on $M$ in that direction into cylinders and minimal components.
Observe that if Flatsurf finds a direction for which $M$ decomposes completely into cylinders, yet those cylinders have incommensurable moduli, then this direction witness that $M$ does not satisfy the Veech dichotomy (See Hubert—Schmidt \cite[Theorem 1]{HubertSchmidt06}). Our code analyzes each of the 92 candidates and in each case finds a direction that is not completely parabolic.
We conclude that none of the 92 candidates generates a Teichm\"uller curve, prooving \autoref{thm:main_thm}.

\bibliographystyle{alpha}
\bibliography{bibli}
\end{document}